\documentclass[a4paper,12pt]{scrartcl}

\usepackage[utf8]{inputenc}

\usepackage{amsmath,amsthm,graphicx,skmath}
\usepackage{bm}
\usepackage[ulem={normalem,normalbf}]{changes}

\usepackage{graphicx,color}

\usepackage{fnpct}

\usepackage{bera}
\usepackage[charter]{mathdesign}

\usepackage{cite}
\usepackage{xspace,microtype}

\DeclareMathOperator{\sign}{sgn}
\newcommand{\diag}{\text{\textbf{diag}}}

\definecolor{darkgreen}{rgb}{0,0.35,0}

\begin{document}

\title{New options for explicit all Mach number schemes by suitable
  choice of time integration methods}

\author{Friedemann Kemm}

\maketitle

\tableofcontents

\begin{abstract}
  Many low-Mach or all-Mach number codes are based on space
  discretizations which in combination with the first order explicit
  Euler method as time integration would lead to an unstable scheme.
  In this paper, we investigate how the choice of a suitable explicit
  time integration method can stabilize these schemes. We restrict
  ourselves to some old prototypical examples in order to find
  directions for further research in this field. 
\end{abstract}

\section{Introduction}
\label{sec:introduction}

One of the most tackled challenges in computational fluid dynamics is
the problem of flows with wide-ranging Mach numbers. While solvers for
low Mach numbers tend to fail for supersonic flows, standard
compressible solvers tend to yield unphysical results for low Mach
numbers. {A first attempt was made by Jameson, Schmidt, and
  Turkel~\cite{doi:10.2514/6.1981-1259} using a special design of the
  numerical viscosity terms. While their attempt already included a
  specially chosen time integration method, it still had flaws, mostly
  when trying to compute steady states.  The next} important step
towards better methods was the asymptotic analysis by Guillard and
Viozat~\cite{viozat} who not only considered the Euler equations of
gas dynamics but also the discretized equations for a first order Roe
scheme. The main finding was that for~\(M\to 0\) the magnitude of the
pressure fluctuations is unphysical, which could only be cured by
lowering the numerical viscosity on acoustic waves.

The first solvers based on these findings were still in the realm of
preconditioning, like the solver by Guillard and
Murrone~\cite{murrone}, but the ideas were transferred over time to
the Riemann solvers themselves as in
\cite{felix-low-mach-fix,philipp-l2roe,CHEN2022111027,doi:10.1137/18M119032X,LI2009810,Yu_2020,HOLMAN2023115129,FLEISCHMANN2020109762}
to name just a few of them.
A comparison of some of these solvers can be found
in~\cite{LI201356}. What many of these solvers have in common is the
reduced numerical viscosity compared to, e.\,g., the Roe solver. As a
consequence, with the explicit Euler method as time discretization,
there is no guarantee for stability\footnote{Most of these papers
  employ higher order space discretizations which considerable reduces
  the influence of the Riemann solver in use. We have already
  discussed this effect in the context of the carbuncle
  phenomenon~\cite{carb15}.}. This has led to a large number of
attempts at implicit or IMEX schemes, i.\,e.\ methods where the
advective terms are discretized explicitly, and the acoustic and, in
the case of the full Navier-Stokes equations, viscous terms are
discretized implicitly. {Also, other ideas of mixing explicit
  and implicit steps can be
  found~\cite{ROSSOW2007879,PELES2018201,Peles2019}.}

A new aspect of these schemes was pointed out quite recently by
Fleischmann et~al~\cite{fleischmann}: In their paper, they were able
to lower the tendency of a code to produce the carbuncle
phenomenon. Since the flow parallel to a shock is at a very low Mach
number, the wrong amplitude of pressure fluctuations in this direction
might trigger the carbuncle. Fleischmann et~el.\ suggested to resort
to a simple version of one of the solvers in~\cite{LI2009810}. In a
previous study~\cite{KEMM2023111947} we discussed this and the
generalized concept of Mach number consistency. In this course we also
introduced a blending between their low-Mach Roe solver and their
simplified model for a traditional carbuncle cure. But even with those
blended solvers no stability is guaranteed when the explicit Euler
method is employed as time integration scheme. {Since in the
  full Navier-Stokes equations the physical viscosity provides a
  stabilizing effect, we restrict our study to the inviscid Euler
  case. It can be expected that stability for Euler flows implies
  stability for Navier-Stokes flows, at least when the time steps are
  small enough for the viscous part.}

The first alternatives to the explicit Euler method were developed in
the mid 19th century by John Couch Adams and tested and published by
Francis Bashforth~\cite{Bashforth2}. Since these methods require
values in the past, they need a starting procedure. To overcome this
difficulty other researchers, especially Carl Runge, Wilhelm Kutta,
and Karl Heun investigated the possibility of higher order one step
methods. The seminal paper of Karl Heun~\cite{Heun1900} marked for
several decades the state of the art in this field. Due to the
contributions by Carl Runge and Wilhelm Kutta these methods are
usually referred to as Runge-Kutta methods. When later the research on
further methods and the stability of the methods in general resurged
again, it was only a few decades until the available material filled
extensive textbooks~\cite{HW1,HW2,HW3,butcherbuch,genlin}. While most
of this research went towards implicit methods, here we want to
restrict ourselves to explicit methods. We want to know which types of
schemes would be best to use in connection with above mentioned low
Mach solvers. A major drawback of explicit Runge-Kutta methods is the
fact that the first stage value is computed via the explicit Euler
scheme. Thus, they might lack the desired stability when applied to
nonlinear problems, whereas multistep methods like Adams-Bashforth do not
suffer from this issue. So, we want to compare some of the oldest
methods of both classes when used in connection with the solver by
Fleischmann et~al.\ or with our blended Mach number consistent
solvers. It can be expected that methods that proof reliable with the
simple Fleischmann solver will also lead to robust schemes when used
in connection with the more elaborate low-Mach solvers mentioned
above. Our goal is to find directions for further research on time
integration methods for low-Mach and all-Mach solvers.

The paper is organized as follows: First we discuss the Riemann solver
by Roe and introduce the simple low-Mach solver by Fleischmann et~al.\
as well as our blended schemes. In Section~\ref{sec:expl-time-integr},
we discuss the basics of explicit Runge-Kutta and Adams-Bashforth
methods with a focus on the linear stability of the schemes and how
this affects the behavior of the fully discretized flow field. After
that we discuss the results for some standard test cases and conclude
the paper with directions for further research in order to optimize
the interplay between space discretization and time integration for
low-Mach and all-Mach number methods. 






\section{The Fleischmann solver and its variants}
\label{sec:fleischm-solv-its}

In this Section, we consider the Riemann solver by Fleischmann et~al.\
as introduced in~\cite{fleischmann} together with our own variants
that blend it with a simple carbuncle fix~\cite{KEMM2023111947}.

\subsection{The Euler equations of gas dynamics}
\label{sec:euler-equations-gas}

The Euler equations for an inviscid gas flow are
\begin{align*}
  \rho_t + \nabla \cdot [\rho \vec v] & = 0\;, \\
  (\rho \vec{v})_t + \nabla \cdot [\rho \vec v \otimes \vec v + p\vec I]  &  =
  0\;, \\ 
  E_t + \nabla \cdot [(E+p)\vec v] & = 0\;
\end{align*}
with the density \(\rho\), the velocity \(\vec v = (u,v,w)^T\) for the
3d-case and \(\vec v = (u,v)^T\) for the 2d-case, the pressure \(p\),
and the total energy \(E\). Density, velocity, and pressure are called
\emph{primitive variables} in contrast to the \emph{conserved
  variables} density, momentum, and total energy. Throughout this
study, we consider the case of an ideal gas with adiabatic coefficient~\(\gamma=1.4\). 

In 2d the flux in \(x\)-direction is
\begin{equation*}
  f(\vec q) = \bigl(\rho u, \rho u^2 + p, \rho u v, (E+p)u\bigr)^T\;.
\end{equation*}
The wave speeds in this direction are the eigenvalues of the according
flux Jacobians: 
\begin{equation}
  \label{eq:2}
  \lambda_1 = u + c\;,\qquad \lambda_2=\lambda_3=u\;,\qquad \lambda_4
  = u + c
\end{equation}
with the speed of sound
\begin{equation*}
  c = \sqrt{\frac{\gamma p}{\rho}} =
  \sqrt{(\gamma-1)\biggl(H-\,\frac{\vec v^2}{2}\biggr)}\;,  
\end{equation*}
where~\(\gamma\) is the ratio of specific heats and
\begin{equation}
  \label{eq:147}
  H = \frac{E+p}{\rho}
\end{equation}
the total enthalpy.

\subsection{Roe-type Riemann solvers}
\label{sec:roe-type-riemann}

In~\cite{roe-orig}, Roe proposed a consistent local linearization
which allows any scalar method to be applied in a wave-wise
manner. Standard upwind on all waves would lead to the well known Roe
flux
\begin{equation}
  \label{eq:4}
  \vec G (\vec q_r,\vec q_l) = \frac{1}{2} \bigl(f(\vec q_r) + f(\vec
  q_l)\bigr) - \frac{1}{2} \abs{\tilde{\vec A}} (\vec
  q_r-\vec q_l)\;, 
\end{equation}
with
\begin{equation*}
  \tilde{\vec A} = \vec A (\vec q_r,\vec q_l)
\end{equation*}
being the consistent local linearization of the flux Jacobian~\(\vec
A\) and the convention
\begin{equation}
  \label{eq:41}
  \abs{\vec A} = \vec R\, \abs{\vec \Lambda}\, \vec L\;, \qquad \abs{\vec
    \Lambda} = \diag\,(\,\abs{\lambda_1},\dots,\abs{\lambda_m}\,)\;. 
\end{equation}
Here~\(\vec L\) and~\(\vec R\) denote matrices of the left and right
eigenvectors respectively. As a result, the~\(\tilde\lambda_k\)
determine the numerical viscosities on the characteristic
fields. These can now be manipulated in order to switch to another
underlying scalar method like HLLE, Rusanov/LLF, etc. A popular
application for such changes are so called entropy fixes, which
prevent the sonic glitch, like the fix by Harten
\begin{equation}
  \label{eq:53}
  \phi(\lambda) =
  \begin{cases}
    \abs{\lambda} & \text{if}~\abs{\lambda} \geq \delta\;,\\
    (\lambda^2+\delta^2)/(2\delta) & \text{if}~\abs{\lambda} <
    \delta\;, 
  \end{cases}
\end{equation}
where \(\delta\) is a small parameter. We will apply this in most
computations of this paper. Only in some cases, we will use the HLLE
speeds~\cite{hllem} on the acoustic waves for the original Roe solver.

\subsection{Mach number consistent Roe solvers}
\label{sec:mach-numb-cons}

First we recall some parts of the discussion we gave
in~\cite{KEMM2023111947}: Guillard and Viozat~\cite{viozat} show that
for very small Mach numbers the viscosity resulting from the wave-wise
upwind as in the standard Roe scheme leads to an incosistency: While
in the low Mach number limit pressure fluctuations scale
with~\(\mathcal O(M^2)\), the numerical scheme supports pressure
fluctuations of order~\(\mathcal O(M)\), where~\(M\) is the reference
Mach number of the flow field. Guillard and Viozat~\cite{viozat}
identify the numerical viscosity as the source of this
inconsistency. To be more precise: They found that the eigenvalues of
the viscosity matrix determine the order of the pressure
fluctuations. Eigenvalues of order~\(\mathcal O(c)\), as we find them
in standard Riemann solvers like Roe on the acoustic waves, lead to
pressure fluctuations of order~\(\mathcal O(M)\). To lower them
to~\(\mathcal O(M^2)\), we need eigenvalues of
order~\(\mathcal O(c\cdot M) = \mathcal O(u)\). This was the starting
point for many of the above mentioned low Mach number and all Mach
number solvers, first of all the solver by Guillard and
Murrone~\cite{murrone}.

A solver that reproduces this behavior can be called Mach number
consistent. This concept can be generalized as follows:
In the low Mach number limit, the viscosities on all waves are of the
same order. This is indeed a generalization of the original concept
since we allow also for higher viscosities as long as they are of the
same order for all waves. Keeping this terminology, Fleischmann et
al.\ provide two basic strategies to maintain Mach number consistency
in the Roe solver: a Mach number dependent upper bound for the
viscosity on the nonlinear, i.\,e.\ acoustic, waves or a lower bound
for the viscosity on the linear waves, i.\,e.\ entropy and shear
waves. With a fixed positive number~\(\phi\), the first approach,
which is the main point of~\cite{fleischmann}, leads to wave speeds
\begin{equation}
  \label{eqfl:1}
  \tilde\lambda_{1,4} = \tilde u \mp \min{\phi\abs{\tilde u},\tilde
    c}\;,\qquad \tilde\lambda_{2,3} = \tilde u\;,
\end{equation}
the second to
\begin{equation}
  \label{eqfl:3}
  \tilde\lambda_{1,4} = \tilde u \mp \tilde c\;,\qquad
  \tilde\lambda_{2,3} =  \sign (\tilde u) \max{\frac{\tilde
    c}{\phi},\abs{\tilde u}}\;. 
\end{equation}
Note that the solver defined by~\eqref{eqfl:1} can be interpreted as a
simple variant of the all-speed Riemann solver by Li, Gu, and
Xu~\cite{LI2009810}, more precisely of their All-Speed-Roe-2.

In~\cite{KEMM2023111947} we suggested two versions of a blended solver
that away from shocks coincide with the one resulting from wave speed
choices~\eqref{eqfl:1}, and in shocks coincide with~\eqref{eqfl:3}.
We suggested the following setting:
\begin{equation}
  \label{eqfl:4}
  \tilde\lambda_{1,4} = \tilde u \mp {\tilde c}^\beta
  \left(\min{\phi\abs{\tilde u},\tilde c}\right)^{1-\beta}\;, \qquad
  \tilde\lambda_{2,3} =  \left( \sign (\tilde u) \max{\frac{\tilde
    c}{\phi},\abs{\tilde u}}\right)^\beta {\tilde u}^{1-\beta}\;,
\end{equation}
which is kind of a weighted geometric mean between the less and the
more viscous approach. Replacing the geometric mean by a weighted
arithmetic mean led us to
\begin{equation}
  \label{eq:8}
  \tilde\lambda_{1,4} = \tilde u \mp \beta \tilde c + (1-\beta)
  \min{\phi\abs{\tilde u},\tilde c}\;, \qquad
 \tilde\lambda_{2,3} = \beta \sign (\tilde u) \max{\frac{\tilde
    c}{\phi},\abs{\tilde u}} + (1-\beta)\tilde u\;.
\end{equation}
Although this would not guarantee full Mach number consistency, our
numerical tests indicate that it is still a reasonable replacement at
lower numerical cost.

\section{Explicit time integration: classical methods}
\label{sec:expl-time-integr}

Here we recall some properties of and some notations for classical
explicit time integration methods. For a more detailed discussion, we
refer to the standard text books~\cite{HW1,HW2,HW3,butcherbuch}.

\subsection{Runge-Kutta methods}
\label{sec:runge-kutta-methods}

A Runge-Kutta method is usually represented by a Butcher-Tableau:
\begin{equation}
  \label{eq:1}
  \begin{array}{c|c}
    \vec c &  \vec A \\
    \hline\\[-.6em]
          & \vec b^T
  \end{array}\quad\text{with}\quad
  \vec A \in \mathbb R^{s\times s},\ \vec b, \vec c \in \mathbb R^s\;,\
  c_i\in[0,1]\ \forall i\;.
\end{equation}
For an ordinary IVP
\begin{equation*}
  \vec y'=\vec f(t,\vec y)\,,\quad \vec f(t_0) = \vec y_o\;,
\end{equation*}
the method can be written as
\begin{align}
  k_i & = f\bigl(t_n+c_j\Delta t, y_n + \Delta t
        \sum_{j=1}^{s} a_{i\,j} k_j\bigr)\quad \text{für}\
        i=1,\dots,s \label{eq:9} \\
  y_{n+1} & = y_n  + \Delta t \sum_{i=1}^{s} b_i k_i \label{eq:6}\\
  \intertext{or as}
  Y_i & = y_n + \Delta t \sum_{j=1}^{s}a_{i\,j}f(tn+c_j\Delta t,
        Y_j) \label{eq:7}\\
  y_{n+1} & = y_n  + \Delta t \sum_{i=1}^{s} b_i f(tn+c_i\Delta t,
            Y_i)\;. \label{eq:10}
\end{align}
For some Runge-Kutta methods, it is also possible to write them as a
kind of sequence of explicit Euler steps. These methods are referred
to as Strong Stability Preserving (SSP), Monotone, Monotonicity
preserving, or sometimes as
TVD~\cite{zbMATH01097356,gottlieb2001strong,zbMATH05820030,zbMATH02139914,zbMATH02220034}. The
advantage is that even for nonlinear ODEs the stability properties of
the explicit Euler scheme are preserved. More precisely, any seminorm
(of the state vector) that would not increase in an explicit Euler
step with the time step size below a certain bound will also not
increase for the monotone scheme with a time step below a constant
multiple of that time step size. The disadvantage is that everything
hinges on the stability of the explicit Euler scheme. As already
mentioned above, this is not the case for the Fleischmann solver and the
blended schemes that we described in Section~\ref{sec:mach-numb-cons}.

In this study, we restrict ourselves to explicit methods, which means
that in the matrix~\(\vec A\) the mean diagonal and the upper triangle
vanish. As can be seen from equations~\eqref{eq:9}--\eqref{eq:10},
this implies that each time step starts with an explicit Euler
step. This rises questions about the stability for nonlinear
problems, especially when the explicit Euler method would not be
stable like for the semidiscrete form of the Euler equations when
discretized by using the Fleischmann solver or some of the other low
Mach solvers cited in the introduction. 

In the following, we further restrict ourselves to methods with
minimal number of stages for the desired order of the scheme. This
means 2nd order two stage, 3rd order three stage, 4th order four stage,
and 5th order 6 stage methods. Since, as will be seen in
Section~\ref{sec:line-stab-meth}, already the 5th order method loses
the main desirable property, we do not consider methods of order six
or higher. 

\subsection{Linear multistep methods}
\label{sec:line-mult-meth}

In general, a linear multistep method has the form
\begin{equation}\label{eq:3}
  \sum_{j=0}^m \alpha_j \vec y_{n-j} = \Delta t  \sum_{j=0}^m \beta_j
  \vec f_{n-j}\qquad\text{with}\quad \vec f_k=\vec f(t_k,\vec y_k)\;. 
\end{equation}
The oldest alternatives to the Euler method are the so called
Adams-Bashforth methods or shorter Adams methods which can be written
in the explicit form
\begin{equation}\label{eq:5}
  \vec y_{n+1} = \vec y_n + \Delta t \sum_{j=0}^s \beta_j \vec f_{n-j}\;. 
\end{equation}
For fixed time step~\(\Delta t\), the coefficients~\(\beta_j\) can be
found in the literature. For variable time steps, the coefficients
also vary. A helpful description can be found
in~\cite[p.~397--400]{HW1}. Note that in order to minimize storage, it
is advantageous to unroll the given recurrence relations.

The concept of Strong Stability Preservation~(SSP) has been extended
also to linear multistep
methods~\cite{zbMATH05015486,zbMATH02192430,zbMATH02027735}. But, as
it still hinges on the stability of the explicit Euler method, it is
not applicable in our case.

A drawback of multistep methods compared to Runge-Kutta methods is the
necessity of a starting procedure since values from the past are
required. In this study, we have avoided this by starting with lower
order schemes---in the first step even explicit Euler---that only rely
on the values already computed or given as the initial value. While
this works for the examples shown below, in general a more elaborate
procedure would be helpful.

\subsection{Linear stability of the methods}
\label{sec:line-stab-meth}

As mentioned in Section~\ref{sec:roe-type-riemann}, wave-wise
application of standard upwind in connection with a Roe-linearization
leads to the standard Roe solver. If we lower the viscosity on the
acoustic waves as in the choice~\eqref{eqfl:1} for the original
Fleischmann solver, we apply to these waves a scalar method with a
lower viscosity than standard upwind. If the flow velocity vanishes,
then the viscosity itself would also vanish which in turn imposes the
application of central differences on the acoustic waves.

\begin{figure}
  \centering
  \includegraphics[width=.4\linewidth]{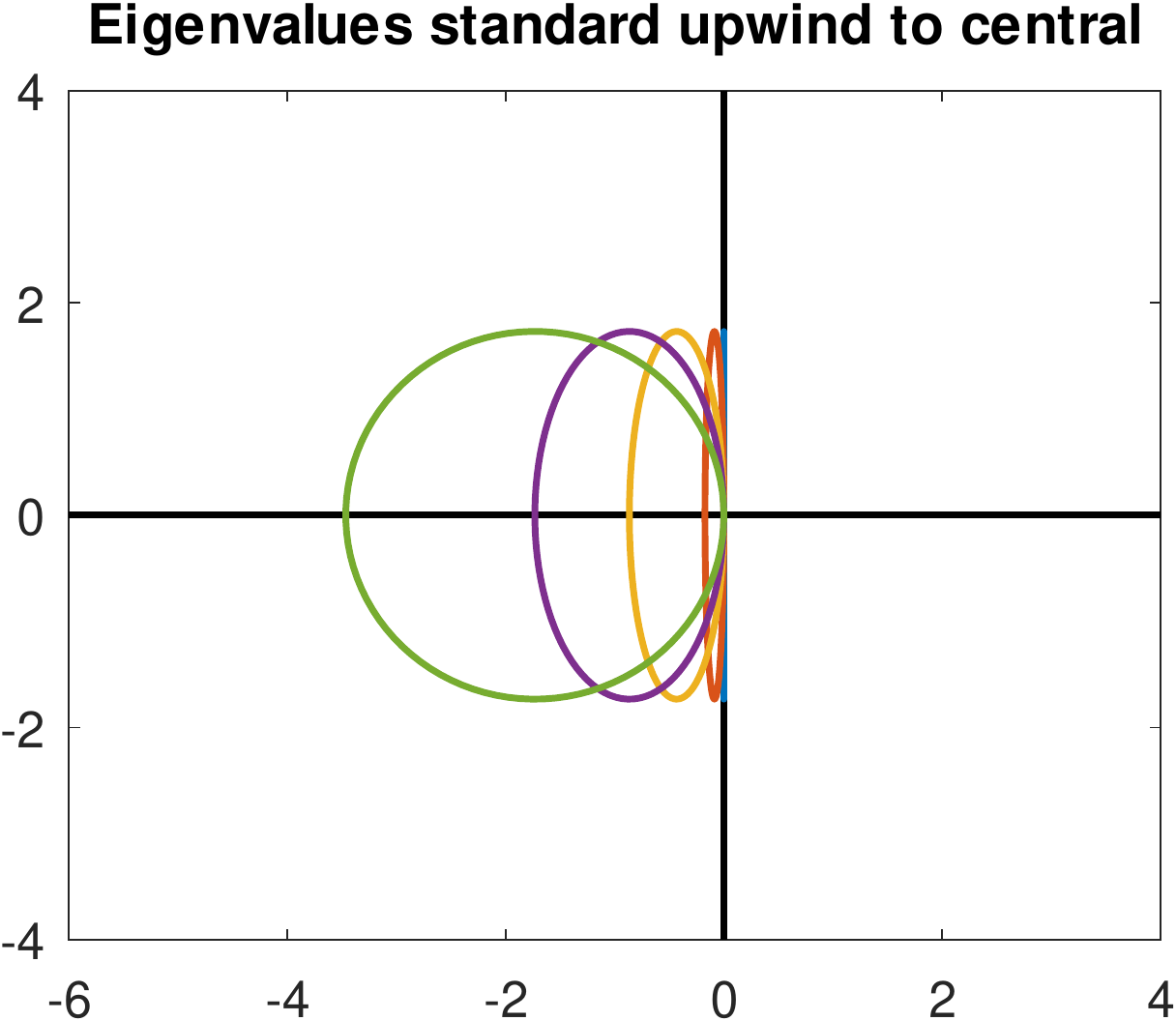}
  \caption[Eigenvalue distribution]{Distribution of eigenvalues for
    space discretization between standard upwind and central.}
  \label{fig:eigdist}
\end{figure}

\begin{figure}
  \centering
  \includegraphics[width=.34\linewidth]{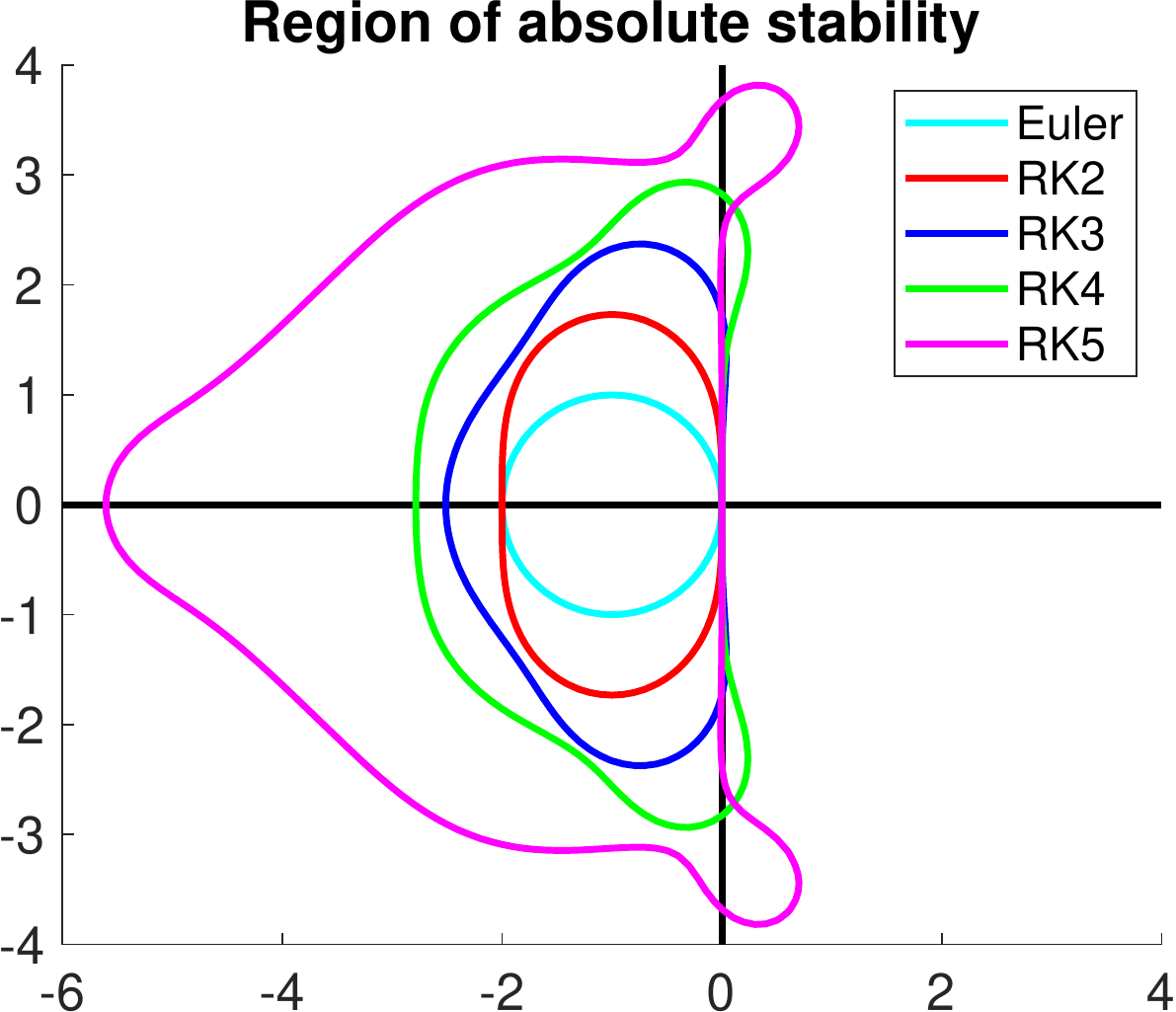}\rule{3em}{0pt}
  \includegraphics[width=.34\linewidth]{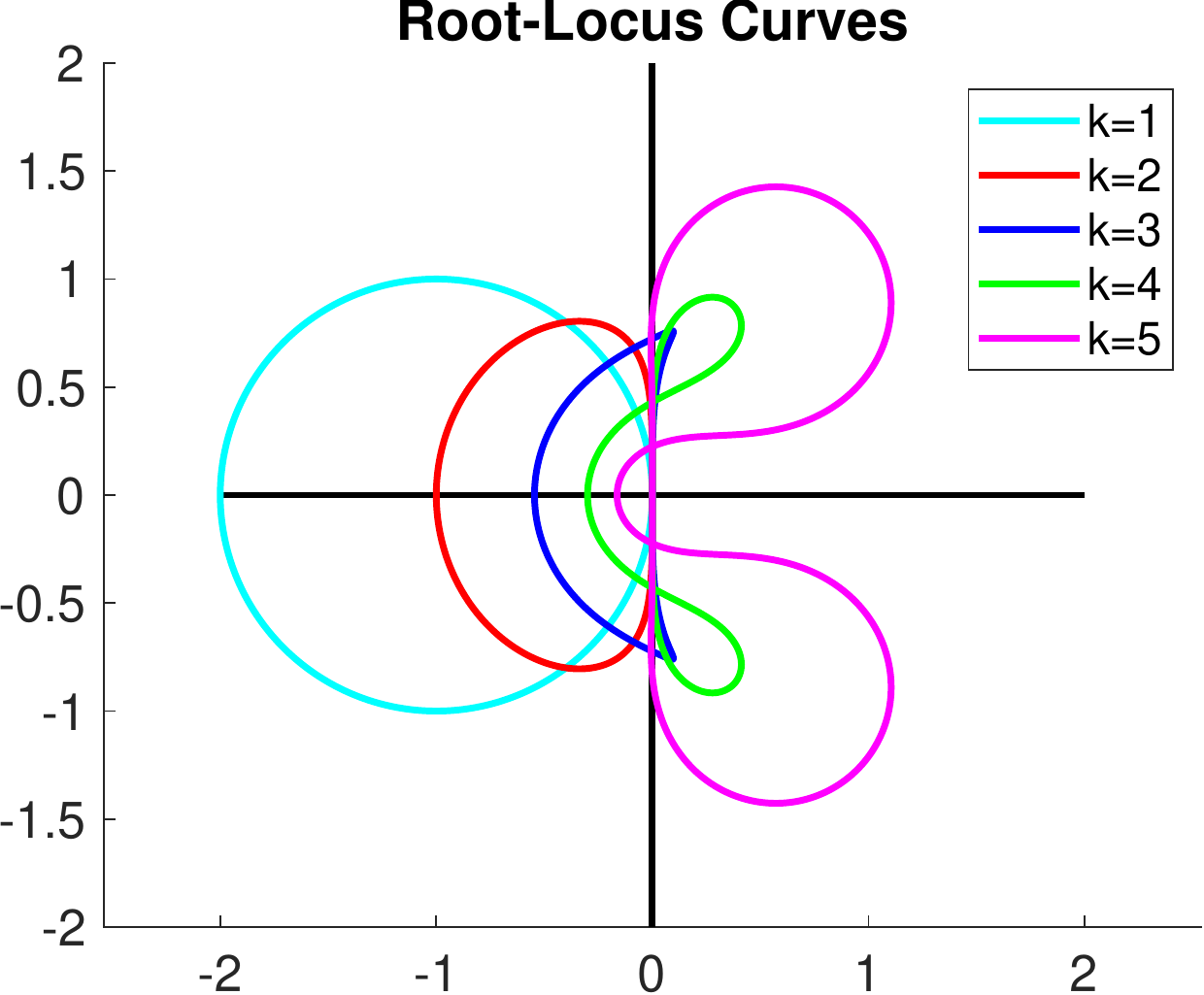}
  \caption[Stability regions]{Stability regions for standard
    schemes. For Adams type schemes, the root locus curves are shown.}
  \label{fig:stbreg}
\end{figure}

In Figure~\ref{fig:eigdist}, we display the eigenvalues of the
discretized linear advection equation using these schemes. While for
standard upwind all eigenvalues are on a circle and for central
differences on a segment of the imaginary axis, for the other schemes
in between, they are located on ellipses. As a consequence, the method
we choose for the time integration has to be stable on the convex hull
of these ellipses and degenerate ellipses. Apparently, the explicit
Euler scheme cannot guarantee stability. Now in
Figure~\ref{fig:stbreg} we have depicted the boundaries of the
stability regions for some standard Runge-Kutta schemes and the root
locus curves for Adams-Bashforth methods with up to five old values
involved. In both plots, the light blue circle marks the stability
region of the explicit Euler method which obviously belongs to both
classes. While in both cases the second order method already promises
better stability, the stability regions do not yet cover the convex
hull of the expected eigenvalues{, except for some
  viscous\footnote{{Viscosity drives the eigenvalues in the
      negative real direction.}}  flows~\cite{grcar2007explicit}.} A
closer look reveals that among the Runge-Kutta schemes only RK3 and
RK4 yield the desired results. Similarly also in the context of Adams
methods, we have to resort to order three or higher.

A consequence of the sizes of the stability regions is that for
Runge-Kutta schemes we can keep the CFL-number from explicit Euler,
even increase it somewhat for RK3 and RK4, but for Adams-Bashforth
methods, we have to reduce the CFL-number with increasing order. Since
for the latter only one flux evaluation is needed per time step, the
computational effort is nearly identical for RK and AB methods of the
same order, at least for the examples we use in this study. Thus the
main remaining question is: How do these schemes perform when applied
to the (nonlinear) Euler equations discretized with the
Fleischmann solver or our blended schemes? The answer will give us
directions for further research and where to look for the ultimate
time integration method.

\section{Numerical tests}
\label{sec:numerical-tests}

Since, as we have seen in the previous section, third order is
sufficient for our purposes, we restrict our comparisons to Euler,
RK2, RK3, and AB2, AB3 (short for Adams-Bashforth of the respective
order). For some tests, we also include the MUSCL-Hancock method,
which obviously is only reasonable in connection with second order in
space. Since for orders two and three (and also four) there are
several two stage and three stage (and four stage) RK-methods
respectively with the same linear stability, we resort to the most
traditional ones, namely those connected to the name of Karl Heun:
\begin{equation*}
  \begin{array}{c|cc}
    0 &   & \\
    1 & 1 & \\
    \hline\\[-.6em]
      & 1/2 & 1/2
  \end{array}\qquad\qquad
  \begin{array}{c|ccc}
    0 &   & & \\
    1/3 & 1/3 & & \\
    2/3 & 0 & 2/3 &\\
    \hline\\[-.6em]
      & 1/4  & 0 & 3/4
  \end{array}    
\end{equation*}
We rechecked with several of the tests and a few other choices for 3rd
order three stage Runge-Kutta but could not find a difference. Thus,
we will stay with only these two methods.

If not stated otherwise, i.\,e.\ for the flow around the infinite
cylinder, in the plots we always show the results for the density.

\subsection{Uniform low-Mach flow}
\label{sec:uniform-low-mach}

As our first example, we consider a grid aligned perturbed low-Mach
flow. The basic flow is a uniform flow in \(x\)-direction but
superimposed with artificial randomized numerical noise. The Mach
number is~\(M=1/20\). The unperturbed initial values for density and
flow velocity (in \(x\)-direction) are~\(\rho = 1,\ u=1\).
The amplitude of the artificial noise on the primitive variables
is~\(10^{-6}\). The final computation time is~\(t=5\).  In the
direction of the flow, we employ first order extrapolation at the
boundaries, in the transverse direction periodic boundary
conditions. Here, we only consider the original Fleischmann solver.
  
Since this test is a simple 2d-extension of a one-dimensional problem,
the results are presented in scatter-type plots: we slice the grid in
\(x\)-direction along the cell faces and plot the density for all
slices at once.

\begin{figure}
  \centering
  \includegraphics[width=.7\linewidth]{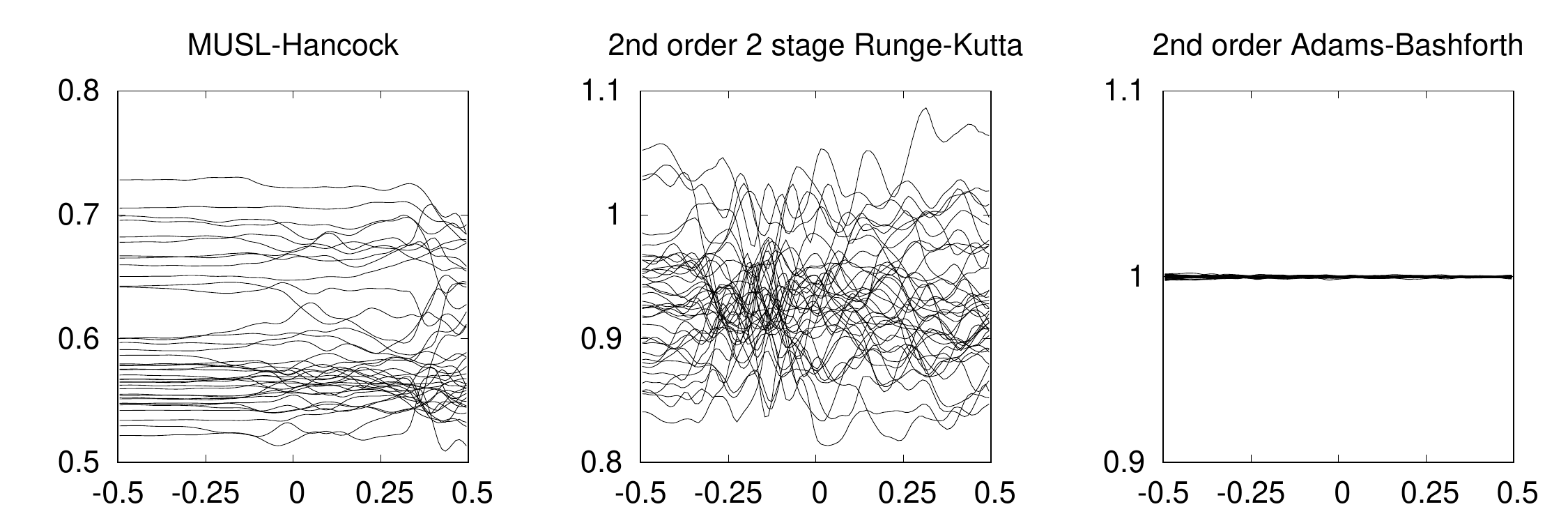}
  \caption[Uniform flow 2nd order]{Perturbed uniform low-Mach flow
    with 2nd order time integration: MUSCL Hancock, 2nd order two
    stage RK (Heun), two step Adams-Bashforth.}
  \label{fig:uniform-o2}
\end{figure}

In Figure~\ref{fig:uniform-o2} results with 2nd order time integration
are shown. It is obvious that AB2 performs superior to MUSCL and
RK2. For the latter two, the perturbations increase and, in connection
with the extrapolation boundaries, even lead to a significant error in
the density.  From linear stability theory, perturbations are expected
to increase with the 2nd order Adams-Bashforth scheme too.  But here
in this case the error is much smaller than for the other two methods.

\begin{figure}
  \centering
  \includegraphics[width=.7\linewidth]{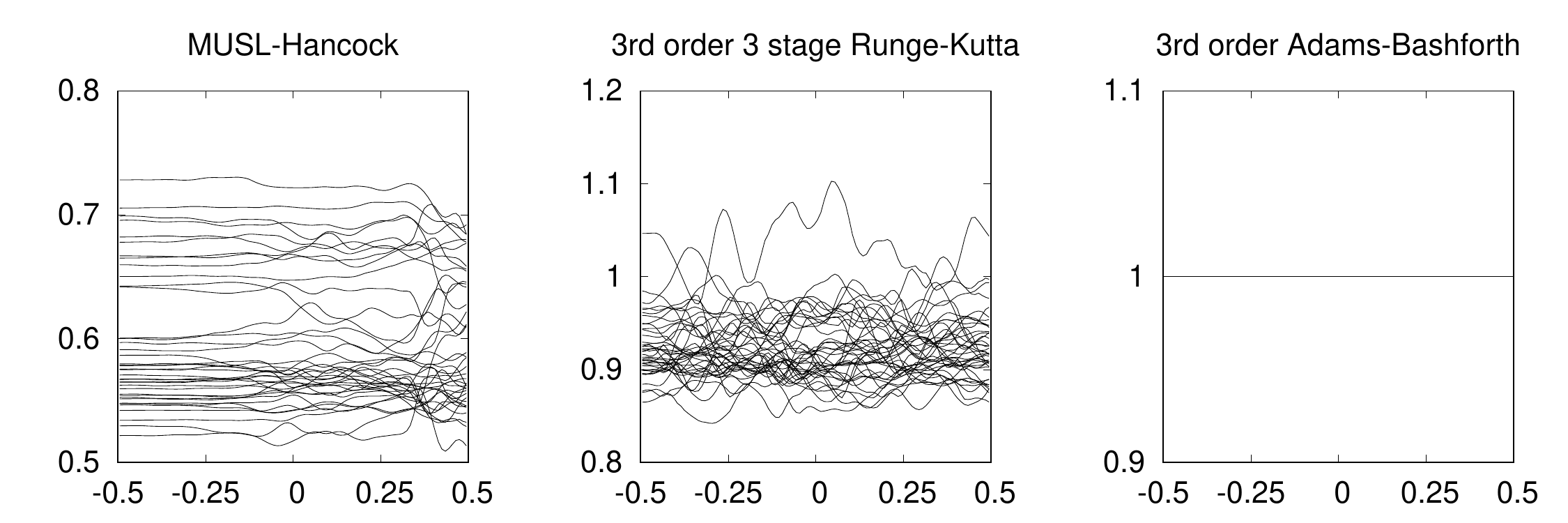}
  \caption[Uniform flow 3rd order]{Perturbed uniform low-Mach flow
    with 3rd order time integration (MUSCL for comparison): MUSCL Hancock, 3rd
    order three stage RK (Heun), three step Adams-Bashforth.}
  \label{fig:uniform-o3}
\end{figure}

Looking at Figure~\ref{fig:uniform-o3}, RK3 shows no significant
improvement compared to RK2. On the other hand, the results with AB3
are perfect as the magnitude of the perturbations does not
increase. This confirms that as we suspected above the starting
explicit Euler step in each time step of an explicit Runge-Kutta
method ruins the stability of the scheme.


\subsection{Flow around infinite cylinder}
\label{sec:flow-around-infinite}

For our second example, we consider the flow around an infinite
cylinder, where we expect the solution for low Mach numbers to
resemble a potential flow with the pressure forming a quadrupole. If
the numerical viscosity on the shear and entropy waves is too high,
this might progressively degenerate, eventually leading to a dipole
structure~\cite{felix-hll-roe}. Another issue is that on quadrilateral
cells the standard Roe solver is not able to capture the physical
solution at all~\cite{felix-et-al-cell}. The reason is that according
to an asymptotic analysis for low Mach numbers, the leading order
velocity shows incompressibility at cell faces which means
that the velocity component perpendicular to the cell face cannot jump
over that face. This should be overcome with the Fleischmann solver
and, hopefully, also with the blended solvers. Again, the issue might
be the stability of the scheme with an unsuited time integration
method.


For the tests, we employ the initial setting~\(\rho = 1, p=1\) with
different Mach numbers. The computational domain is an annulus with
inner radius~\(r=1\) and outer radius~\(R=5\) discretized with
160~cells in circumferential and 100~cells in radial direction. At the
cylinder we assume wall conditions, at the outer boundary Dirichlet
conditions set to the initial state. With these settings, the solution
eventually becomes stationary. In the following, the plots show the
pressure field for the stationary discrete flow. All computations are
done with first order in space.

\begin{figure}
  \centering
  \includegraphics[width=.9\linewidth]{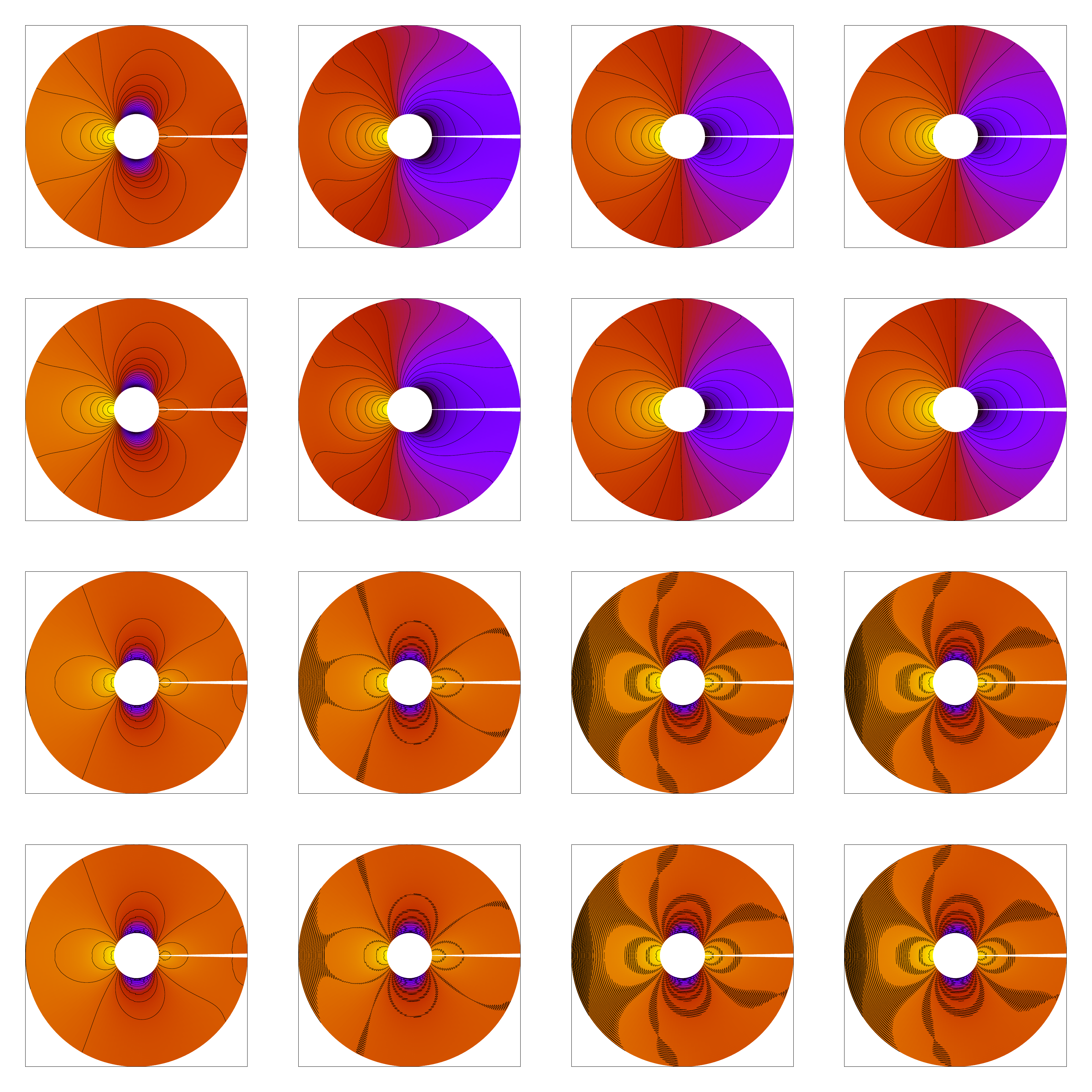}
  \caption[Cylinder flow AB\,3]{{Subsonic flow around infinite cylinder
    with 3rd order Adams scheme for Mach
    numbers~\(10^{-1},\,10^{-2},\,10^{-3},\) and~\(10^{-4}\) (left to
    right) with original Roe and 3rd order AB (first row) and RK
    (second row) as well as the low-Mach solver suggested by
    Fleischmann et~al.\ also with 3rd order AB (third row) and 3rd
    order RK (last row).}} 
  \label{fig:lm-all-ab3}
\end{figure}

We start with the 3rd order time integration
methods. {Figure~\ref{fig:lm-all-ab3} shows results with the
  original Roe solver and the low-Mach solver by Fleischmann
  et~al. While the results with the original Roe solver are unphysical
  as expected for very low Mach numbers, the results with the low-Mach
  solver are reasonable. Only slight fluctuations can be seen that
  would be invisible if they were not highlighted by the contour
  lines. Interestingly, for this example both time integration methods
  behave similar. This might explain why this issue, which is
  introduced by the unstable first stage value of explicit Runge-Kutta
  schemes, was not discovered earlier.}

\begin{figure}
  \centering
  \includegraphics[width=.4\linewidth]{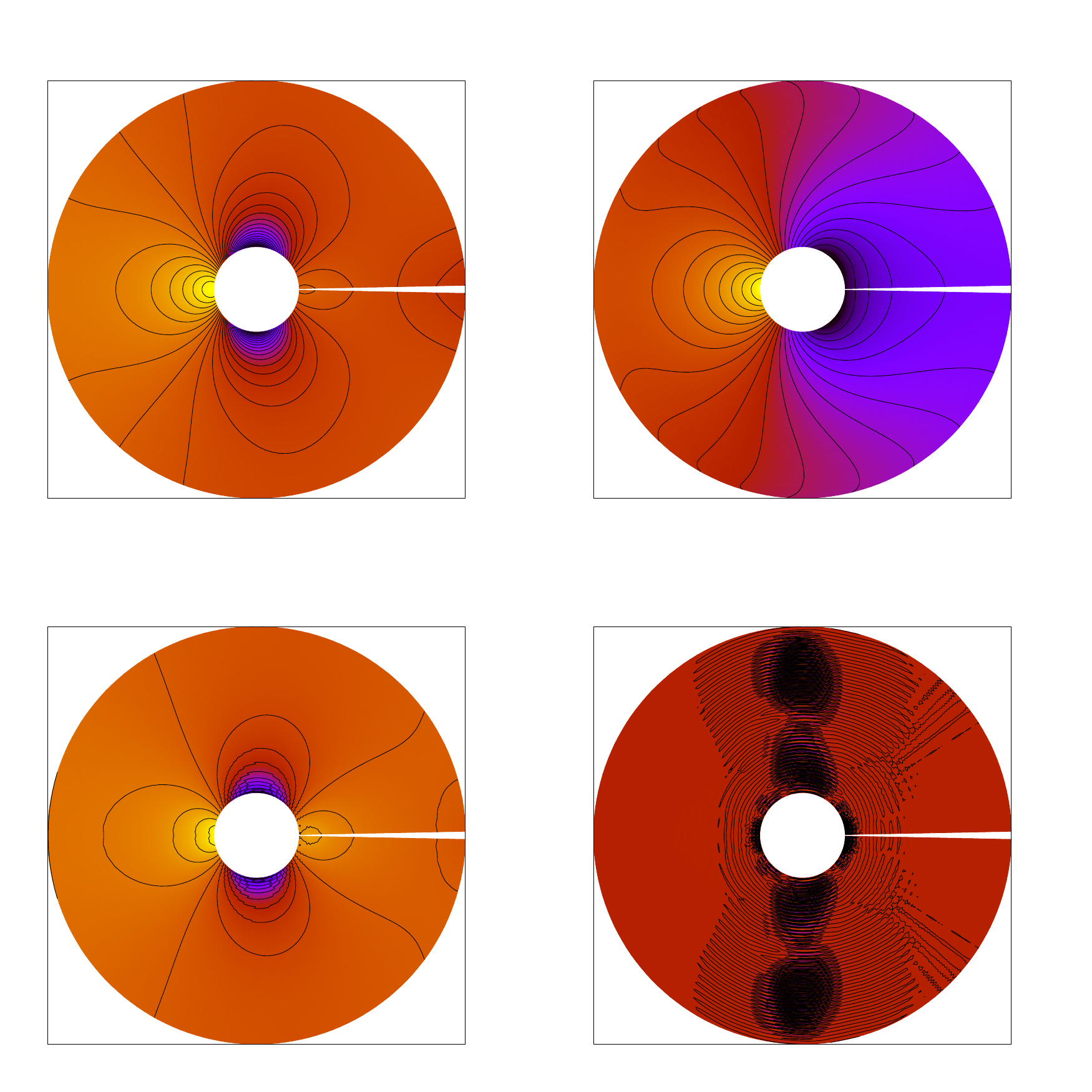}
  \caption[Cylinder flow RK\,2]{Subsonic flow around infinite cylinder
    with 2nd order Runge-Kutta scheme for numbers~\(10^{-1}\) (left)
    and~\(10^{-2}\) (right) with original Roe solver (top row) and
    Fleischmann type solvers (bottom row).}
  \label{fig:lm-all-rk2}
\end{figure}

\begin{figure}
  \centering
  \includegraphics[width=.45\linewidth]{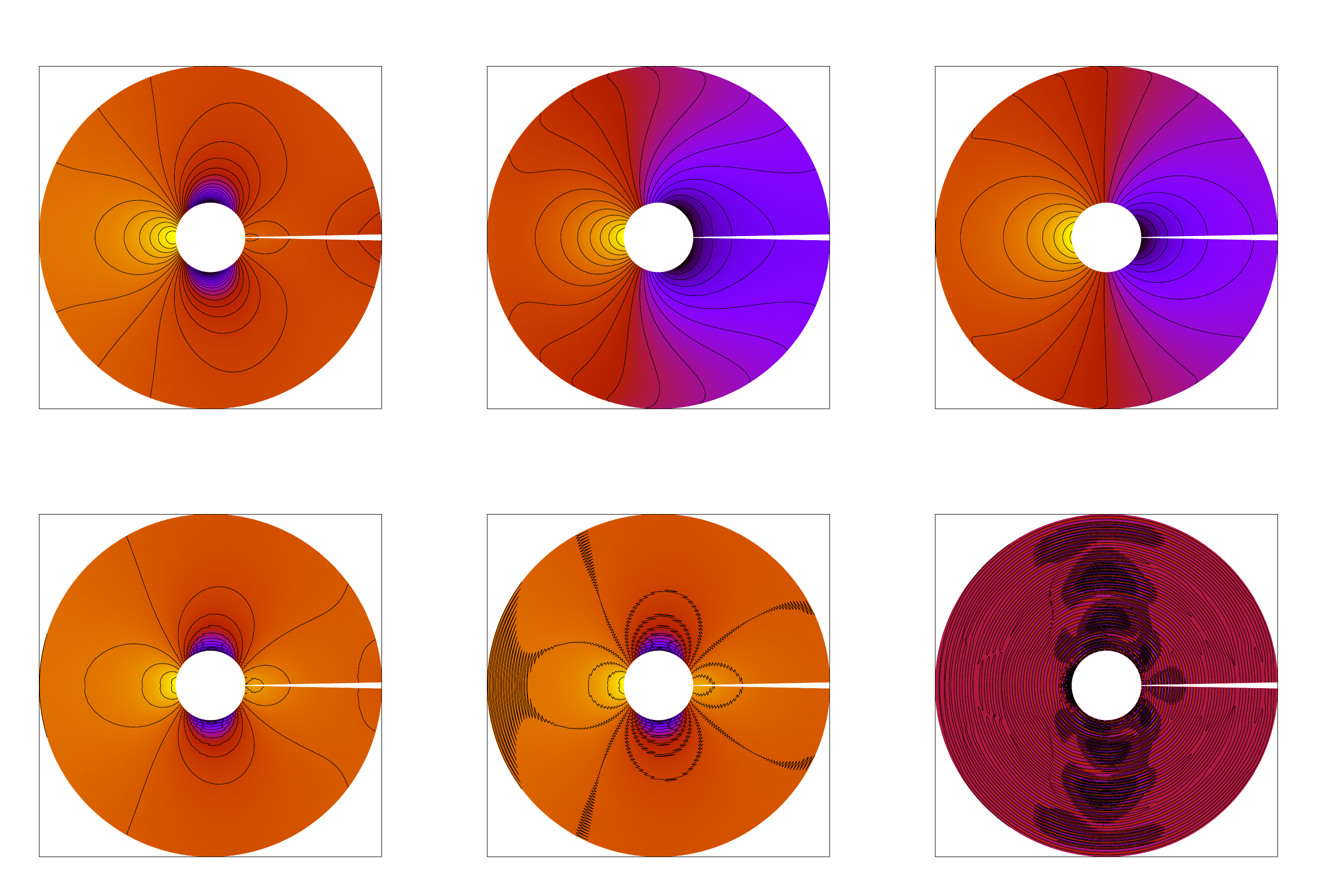}
  \caption[Cylinder flow AB\,2]{Subsonic flow around infinite cylinder
    with 2nd order Adams scheme for Mach
    numbers~\(10^{-1},\,10^{-2},\) and~\(10^{-3}\) (left to
    right) with original Roe (top row) and
    Fleischmann type solvers (bottom row).}
  \label{fig:lm-all-ab2}
\end{figure}

The results with the second order methods are shown in
Figures~\ref{fig:lm-all-rk2} and~\ref{fig:lm-all-ab2} for Runge-Kutta
and Adams-Bashforth respectively. {Due to the unstable first
  stage value, the results for the Runge-Kutta method are inferior to
  those obtained with the linear multistep method. But, as expected
  from the stability region, they are also inferior to the results
  with 3rd order Adams-Bashforth.} Already for~\(M=10^{-3}\) the
solution is completely worthless due to heavy oscillations.

\subsection{Supersonic blunt body flow}
\label{sec:supers-blunt-body}

Since we are interested in all-Mach solvers, we also have to include
supersonic problems. As a first example we consider a~\(M=20\)
blunt-body flow. As in our previous work~\cite{KEMM2023111947}, we
choose an infinite cylinder as test situation, more precisely, a
cylinder with radius~\(r=1\). As computational domain we employ a
third of the annulus with inner radius~\(r=1\) and outer
radius~\(R=2\). Since the interesting part of the flow is the inflow
region, we restrict the domain in angular direction
to~\([{2\pi/3},{4\pi/3}]\). The domain is
discretized with 150 cells in the radial direction and 800 cells in
the angular direction. The initial flow is set to the inflow state
everywhere. At the cylinder, we employ wall boundary conditions, at
the other boundaries first order extrapolation.

\begin{figure}
  \centering
  \includegraphics[width=.7\linewidth]{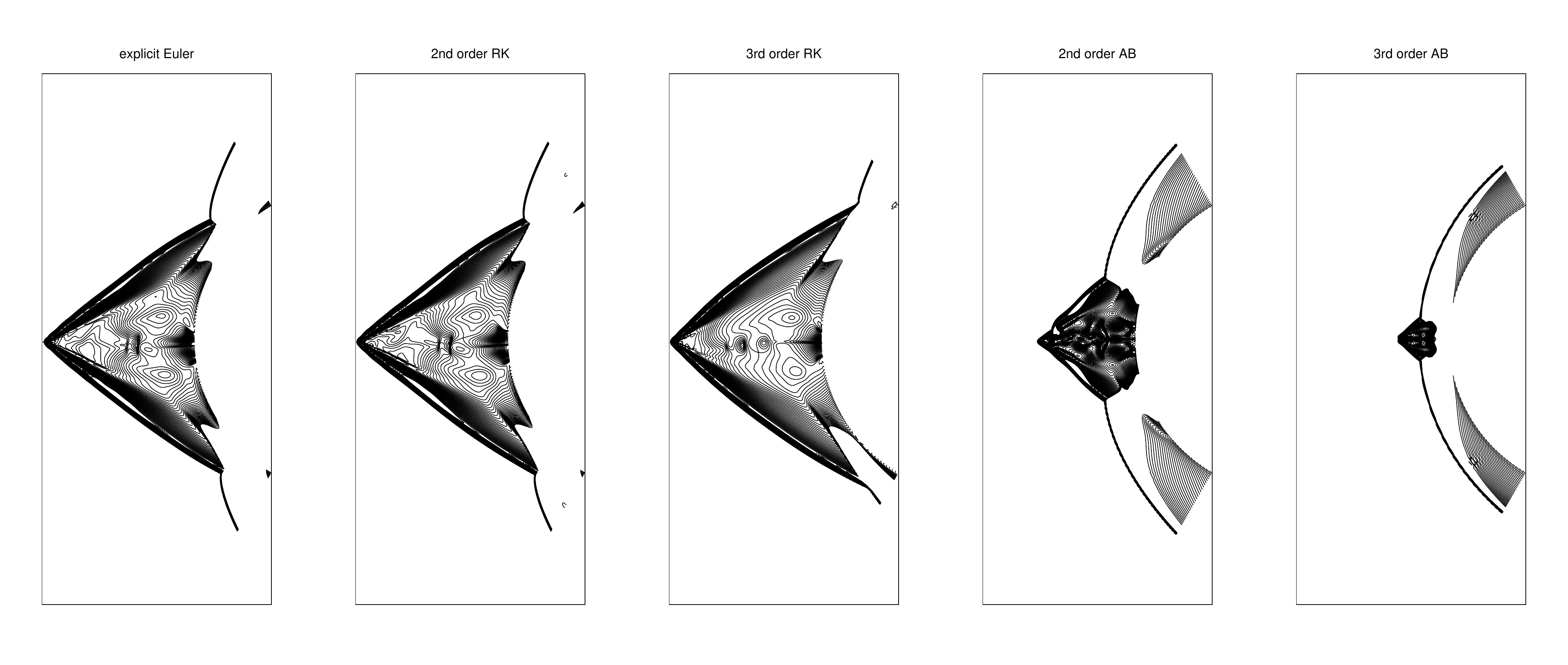}
  \caption[Blunt body flow different time schemes]{Mach 20 blunt body
    flow at time~\(t=0.5\) with Fleischmann solver and different time
    integrations: Explicit Euler, Runge-Kutta 2nd and 3rd order,
    Adams-Bashforth 2nd and third order.}
  \label{fig:eerkab}
\end{figure}

\begin{figure}
  \centering
  \includegraphics[width=.55\linewidth]{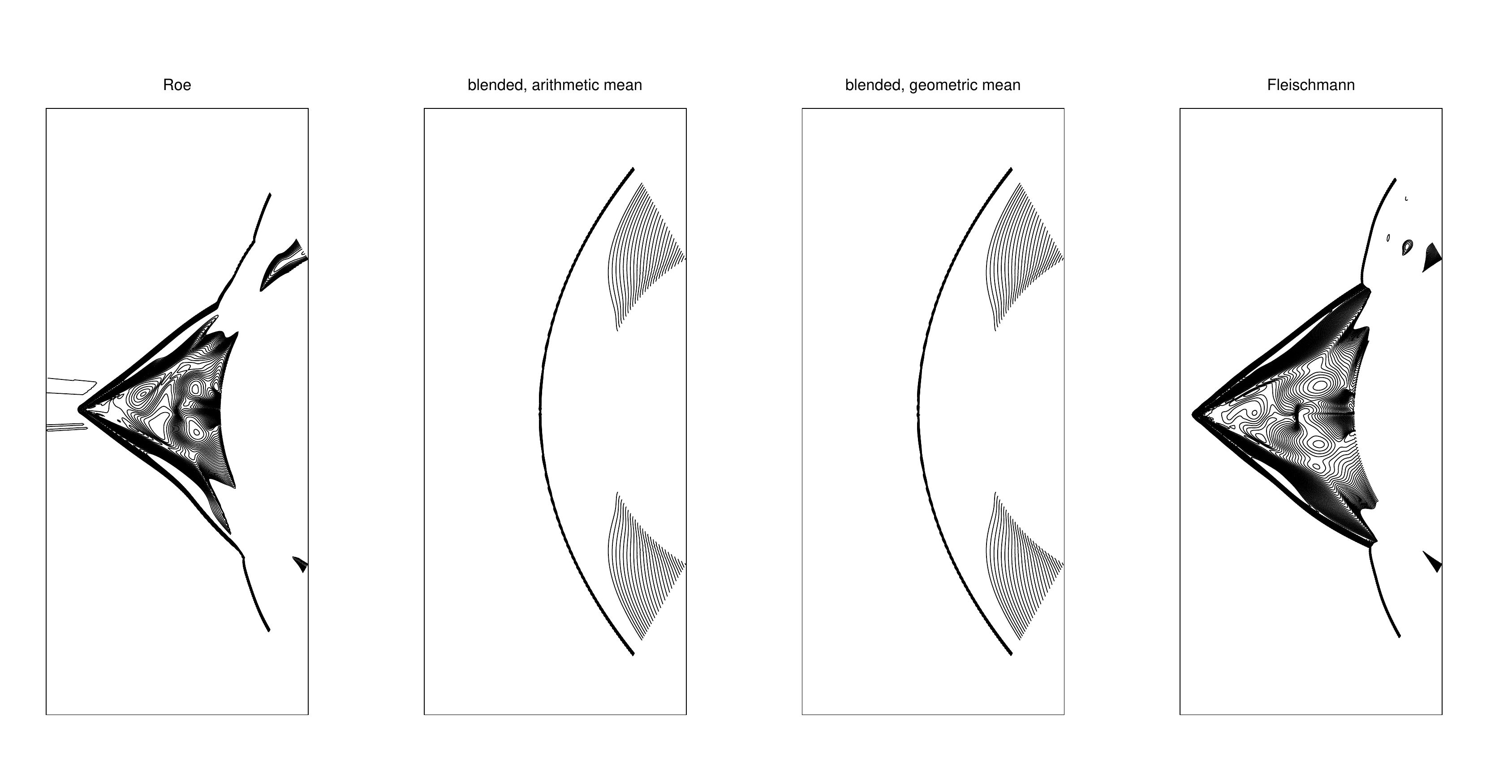}
  \caption[blunt body flow AB\,3]{Mach 20 blunt body flow at time~\(t=2\)
    with 3rd order Adams-Bashforth and different Riemann solvers: Roe,
    blended schemes with arithmetic and geometric mean, Fleischmann
    solver.}
  \label{fig:fl-ab3}
\end{figure}

The original purpose of Fleischmann et~al.~\cite{fleischmann} was to
prevent the carbuncle phenomenon by reducing perturbations in the
transverse direction, where the directional Mach number is small. In
our case, instabilities of the scheme may also lead to
perturbations. Thus, we first compare the different time integration
methods in connection with the original Fleischmann solver at an early
time. The results are displayed in Figure~\ref{fig:eerkab}. Obviously,
the explicit Euler method and the explicit Runge-Kutta methods, which
include an explicit Euler step for the computation of the first stage
value, cannot prevent the carbuncle. The results for these schemes are
nearly identical. Only with the linear multistep methods, the
situation improves. But even with the 3rd order Adams-Bashforth method
a carbuncle will form.

In Figure~\ref{fig:fl-ab3}, we show the results for a later time with
AB3 and different Riemann solvers: original Roe, the blended solvers,
and the original Fleischmann solver. Obviously, as was expected
from~\cite{KEMM2023111947}, only the blended solvers actually prevent
the carbuncle. This means that even with AB3 perturbations are still
created in the transverse direction which trigger the instability of
the discrete shock profile.

\subsection{{Richtmyer-Meshkov instability}}
\label{sec:richtmy-meshk-inst}

\begin{figure}
  \centering
  \includegraphics[width=.5\linewidth]{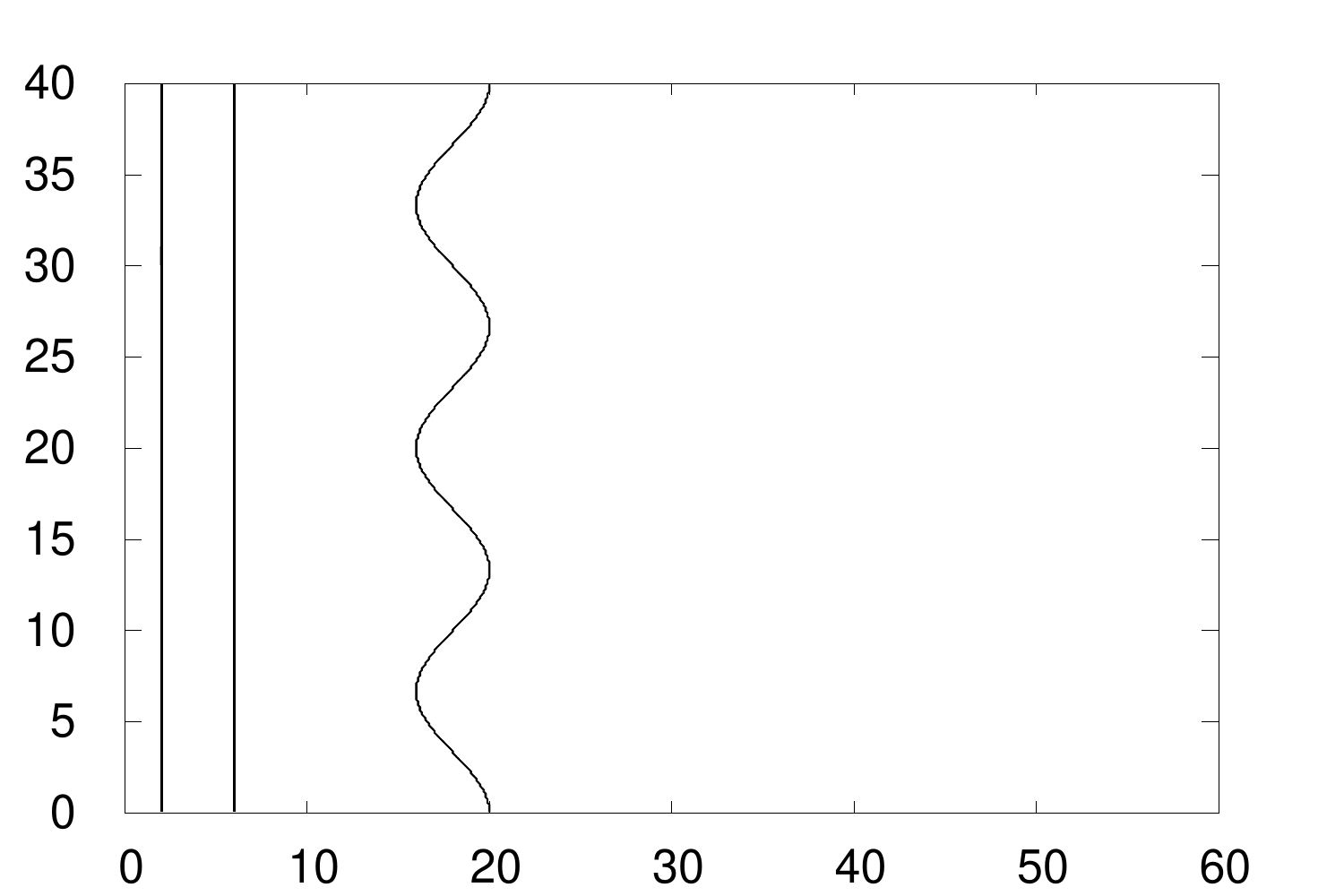}
  \caption{{Initial setting for RMI: In strip \(x\in [2,6]\) higher
    pressure and density, namely~\(p=4.9,\) and~\(\rho=4.22\).
    Otherwise~\(p=1\), density left and right of sinusoidal
    boundary~\(\rho=1\) and~\(\rho=0.25\).}}
  \label{fig:iniRMI}
\end{figure}

\begin{figure}
  \centering
  \includegraphics[width=.75\linewidth]{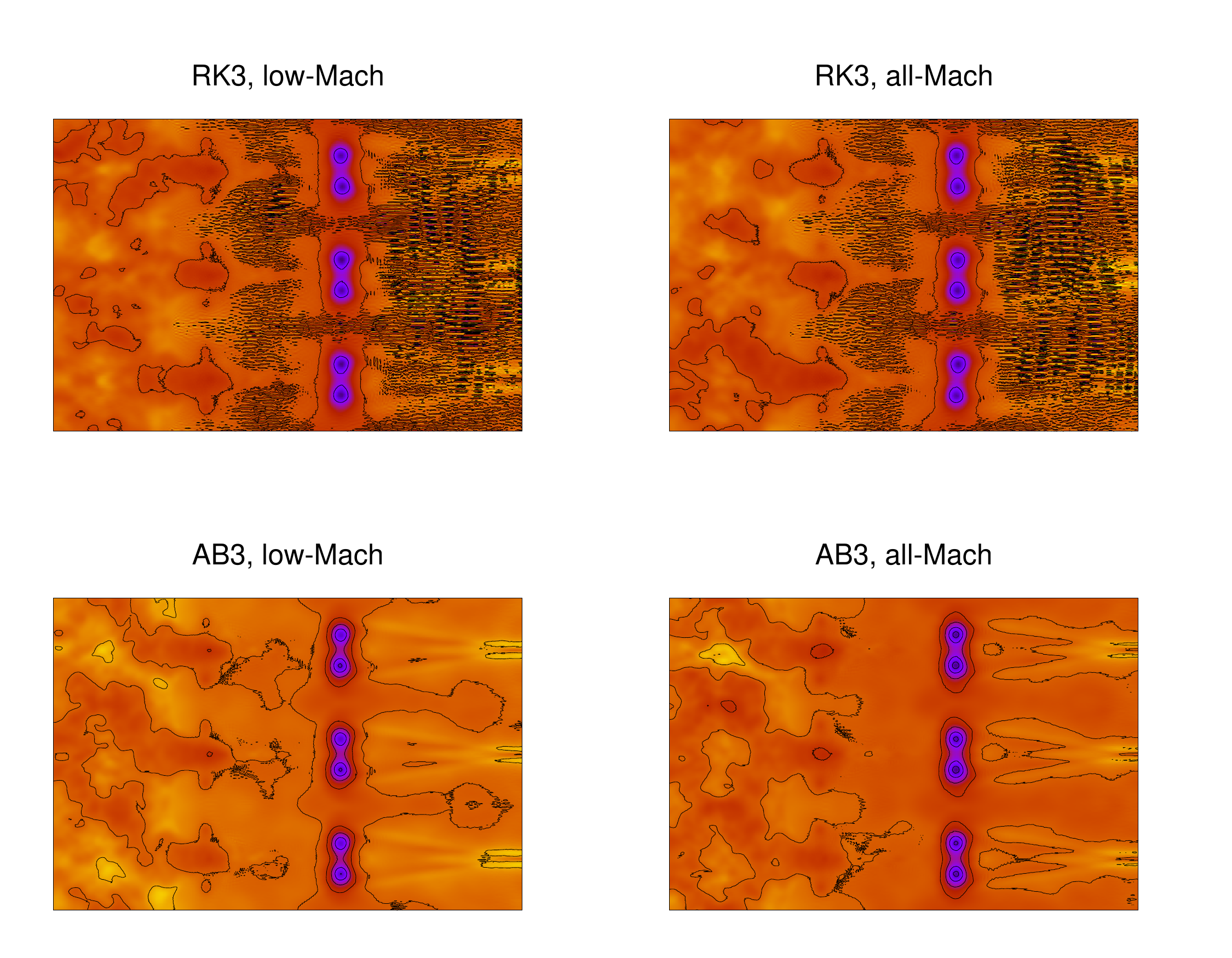}
  \caption{{Richtmyer-Meshkov instability, pressure at time~\(t=85\)}}
  \label{fig:rmi}
\end{figure}

The so called Richtmyer-Meshkov
instability~\cite{https://doi.org/10.1002/cpa.3160130207} results in a
flow field that is comprised of different regimes, clearly
compressible regions and nearly incompressible parts. Initially, the
flow is at rest. Pressure and density are as sketched out in
Figure~\ref{fig:iniRMI}. This causes a shock that runs over the
sinusoidal phase boundary and, thus, drives the dynamics of the flow
field that eventually leads to a mixing of both the lighter and the
heavier fluids. In Figure~\ref{fig:rmi}, we show numerical results for
the pressure with low-Mach solver as suggested by Fleischmann (left
column) and the blended Mach-number consistent solver with the
geometrical mean as an all-Mach solver (right column) with both third
order 3-stage Runge-Kutta (upper row) and 3rd order Adams-Bashforth
(lower row). All computations are with second order in space and a
resolution of 600 times 396 cells. The differences between the two
Riemann solvers are small. The differences between the time
integration methods are, however, much more prominent. Expecially in
the low density part with an extremely low Mach number, the results
with Runge-Kutta show significant instabilities, as expected from the
instability of the first stage value and the results from
Section~\ref{sec:uniform-low-mach}.%

\subsection{Test examples that leave open questions: Elling and Kelvin-Helmholtz}
\label{sec:failure-some-test}

\begin{figure}
  \centering
  \includegraphics[width=.45\linewidth]{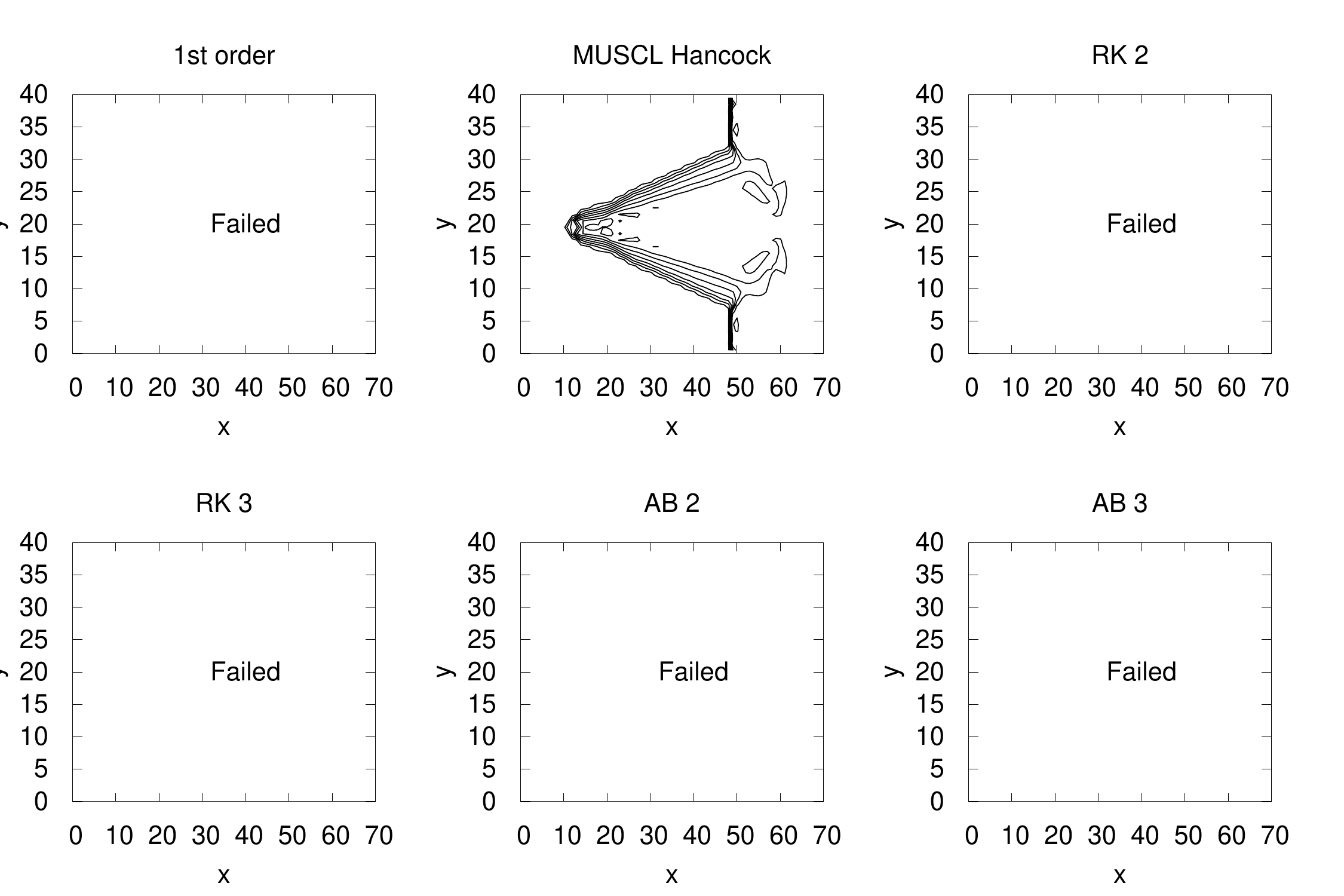}
  \caption[Elling test]{Elling test with Fleischmann solver and
    different time integration methods. Except MUSL-Hancock, all fail.}
  \label{fig:ell}
\end{figure}

While the previous results clearly favor Adams methods, there are some
examples where no difference exists between Runge-Kutta and
Adams-Bashforth like the Elling test, which is a test for a solvers
ability to reproduce physical
carbuncles~\cite{elling_carbuncle_2009}. As is
indicated in Figure~\ref{fig:ell}, while the original Fleischmann
solver performs nicely with MUSCL-Hancock, neither the Runge-Kutta
methods nor the Adams-Bashforth methods provide enough stabilization
to avoid failure of the computation due to negative internal
energy. Even the blended solvers fail. 

\begin{figure}
  \centering
  \includegraphics[width=.75\linewidth]{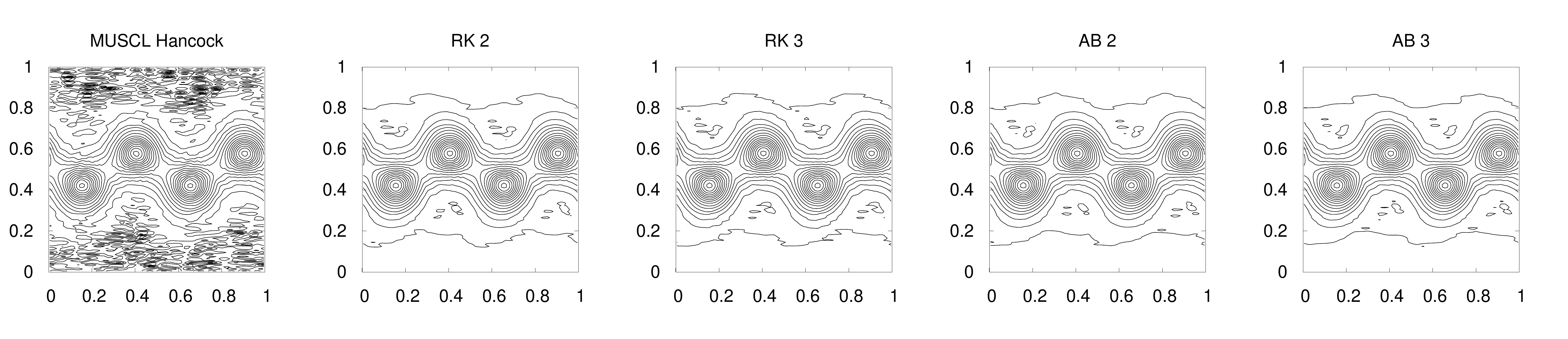}
  \caption[Kelvin-Helmholtz]{Density for Kelvin-Helmholtz instability
    with Fleischmann solver and different time integration methods at
    time~\(t=4\), contour plot}
  \label{fig:KHIcontour}
\end{figure}

On the other hand, the computation of a Kelvin-Helmholtz
instability (Figure~\ref{fig:KHIcontour}) yields results for Runge-Kutta
and Adams-Bashforth that are correct and nearly indistinguishable while
MUSCL-Hancock suffers from the expected instability due to the low
viscosity of the Fleischmann solver.

\section{Conclusions and possible directions for further research}
\label{sec:concl-direct-furth}

We tested some classical explicit time integration methods with the
simple low Mach number Riemann solver by Fleischmann et~al.\ as well
as with its blended versions. It became apparent that explicit
Runge-Kutta methods are not the best choice: The fact that the first
stage value is computed with explicit Euler {spoils} the
stability of the resulting scheme. As an alternative we tested some
simple Adams-Bashforth methods which behaved considerably better. For
many of the more elaborate low-Mach and all-Mach solvers cited in the
introduction, they might be sufficient, especially with higher order
in space. With very high orders even Runge-Kutta methods may do the
job. But for the general case we have to look for further
improvement. Some possibilities would be PECE-type methods or general
linear methods~\cite{genlin}. For the latter, the stability of stage
values has to be taken into account since unstable stage values still
could ruin the stability of the scheme.


\bibliographystyle{amsplain} \bibliography{allmach-cc}

\end{document}